\documentclass[12pt,reqno]{amsart}
\usepackage{amssymb,delarray}
\usepackage{amsfonts,amscd}
\usepackage{epsfig}
\usepackage[all]{xy}

\def\geq{\geqslant}

\newtheorem{thm}{Theorem}
\newtheorem{prop}[thm]
{Proposition}

\newtheorem{ex-thm}[thm]{Theorem-Example}

{\catcode`\@=11
\gdef\n@te#1#2{\leavevmode\vadjust{%
 {\setbox\z@\hbox to\z@{\strut#1}%
  \setbox\z@\hbox{\raise\dp\strutbox\box\z@}\ht\z@=\z@\dp\z@=\z@%
  #2\box\z@}}}
\gdef\leftnote#1{\n@te{\hss#1\quad}{}}
\gdef\rightnote#1{\n@te{\quad\kern-\leftskip#1\hss}{\moveright\hsize}}
\gdef\?{\FN@\qumark}
\gdef\qumark{\ifx\next"\DN@"##1"{\leftnote{\rm##1}}\else
 \DN@{\leftnote{\rm??}}\fi{\rm??}\next@}}

\begin{document}

\title{%
The Hesse curve of a Lefschtz pencil of plane curves}

\author{Vik.\,S.~Kulikov}
\date{}

\address{Steklov Mathematical Institute}
 \email{kulikov@mi.ras.ru}

\begin{abstract} 
We prove that for a generic Lefschetz pencil of plane curves of degree
$d\geq 3$ there exists a curve $H$ (called the  Hesse curve of the pencil) of degree
$6(d-1)$ and  genus $3(4d^2-13d+8)+1$, and such that: $(i)$ $H$ has $d^2$ singular points of multiplicity three at the base points of the pencil
and $3(d-1)^2$ ordinary nodes at the singular points of the degenerate members of the pencil;
$(ii)$ 
 for each member of the pencil the intersection of $H$ with this fibre consists of the inflection points of this member and the base points of the pencil.
\end{abstract}

\maketitle

\setcounter{tocdepth}{1}

\def\st{{\sf st}}

\setcounter{section}{-1}

{\bf 1.} Let  $F(\overline a,\overline z)=\displaystyle \sum_{k+m+n=d}a_{k,m,n}z_1^kz_2^mz_3^n$ be the homogeneous polynomial of degree $d$ in variables $z_1,z_2,z_3$ and of degree one in variables $a_{k,m,n}$, $k+m+n=d$. Denote by $\mathcal C_d\subset \mathbb P^{K_d}\times \mathbb P^2$, where $K_d=\frac{d(d+3)}{2}$, the complete family of plane curves of degree $d$ given by equation $F(\overline a,\overline z)=0$.
     Let ${\mathcal I_d}=\mathcal C_d\cap\mathcal H_d$, where
$$\displaystyle \mathcal H_d=\{ (\overline a,\overline z)\in \mathbb P^{K_d}\times \mathbb P^2\mid\det (\frac{\partial^2 F(\overline a,\overline z)}{\partial z_i\partial z_j})=0\} .$$

Denote by $f_d:\mathcal C_d\to  \mathbb P^{K_d}$ the restrictions of the projection $\text{pr}_1: \mathbb P^{K_d}\times \mathbb P^2\to \mathbb P^{K_d}$ to $\mathcal C_d$. It is well-known (see, for example, \cite{BK}) that for a generic point $\overline a_0\in \mathbb P^{K_d}$ the intersection of the curve $C_{\overline a_0}=f_d^{-1}(\overline a)$ and its Hessian curve $H_{C_{\overline a_0}}$ given by $\frac{\partial^2 F(\overline a_0,\overline z)}{\partial z_i\partial z_j})=0$ is the set of the inflection points of $C_{\overline a_0}$ containing $3d(d-2)$ points. Therefore  for $d\geq 3$  the morphism $h_d= f_{d|\mathcal I_d}:{\mathcal I_d}\to \mathbb P^{K_d}$  has degree $\deg h_d=3d(d-2)$.

Let $\mathcal S_d$ be a subvariety of $\mathbb P^{K_d}$ consisting of the points $\overline a$ such that the curves $C_{\overline a}$ are singular and let $\mathcal M_d$ be a subvariety of $\mathbb P^{K_d}$ consisting of the points $\overline a$ such that for $\overline a\in \mathcal M_d$ the curve $C_{\overline a}$ has a $r$-tuple inflection point with $r\geq 2$.
Let $\mathcal B_d=\mathcal S_d\cup \mathcal M_d$ (if $d=3$ then $\mathcal M_3=\emptyset$). It is easy to show (\cite{Ku4}) that $\mathcal M_d$ is an irreducible hypersurface in $\mathbb P^{K_d}$ if $d\geq 4$.
It is well-known also that $\mathcal S_d$ is an irreducible hypersurface in $\mathbb P^{K_d}$,
$\deg \mathcal S_d=3(d-1)^2$.
\begin{prop} {\rm (\cite{Ku4})} \label{propo}
The local monodromy group \footnote{The definition of the local monodromy group of a dominant morphism $\varphi: X\to \mathbb P^m$ at a point $p\in\mathbb P^m$ can be found, for example, in \cite{Ku2} or \cite{Ku3}.} of $h_d$ at a generic point $\overline a\in \mathcal M_d$ is a subgroup $\mathbb Z_2$ of the symmetric group $\mathbb S_{3d(d-2)}$ generated by a transposition, and the local monodromy group at a generic point $\overline a\in \mathcal S_d$ is a subgroup $\mathbb Z_3\subset\mathbb S_{3d(d-2)}$ generated by the product of two disjoint cycles of length three.
\end{prop}

{\bf 2.} Remind that, by definition, a {\it Lefschetz pencil} is a fibration $f_d: C_L=f_d^{-1}(L)\to L$
of curves over $L$ (and also the linear system $\{ C_{\overline a}\}_{\overline a\in L}$ of curves of degree $d$ in $\mathbb P^2$), where $L$ is a line in $\mathbb P^{K_d}$ in general position with respect to the divisor $\mathcal S_d$. The  body $C_L$ of the Lefschetz pencil is a non-singular surface. The linear system $\{ C_{\overline a}\}_{\overline a\in L}$ has  $d^2$ base points and the restriction of the projection $pr_2$ to $C_L$ is the composition of $d^2$ $\sigma$-processes with centers at the base points.
We say that the Lefschtz pencil $f_d: C_L\to L$ is {\it generic} if $L$ and $\mathcal M_d\setminus \mathcal S_d$ meet at $m_d=\deg \mathcal M_d$ different points.

\begin{prop} \label{m_n} We have $\deg \mathcal M_d=6(d-3)(3d-2)$.
\end{prop}
\proof Consider a generic Lefschetz pencil $\{ C_{\overline a}\}_{\overline a\in L}$. It follows from Theorem 1 in \cite{Ku4} that the curve $I_L=h_d^{-1}(L)=C_L\cap H_{C_L}$ is irreducible.
Also, it is easy to show (\cite{Ku4}) that the curve $I_L$ has $\deg S_d=3(d-1)^2$ singular points which are the ordinary nodes. Let
$\nu :\widetilde I_L\to I_L$ be the normalization of $I_L$, $g_d$ the genus of $\widetilde I_L$, and $\widetilde h_d=h_d\circ \nu: \widetilde I_L\to L$. It follows from Proposition \ref{propo} that $\widetilde h_d:\widetilde I_L\to L$ is ramified at $m_d$ points with multiplicity two (at the preimages of $2$-tuple inflection points of the fibres of the Lefschetz pencil) and it is ramified at $2\deg \mathcal S_d=6(d-1)^2$ points with multiplicity three (at the preimages of the singular points of $I_L)$. Therefore, it follows from Hurwitz formula that
\begin{equation}\label{g1} 2(g_d-1)=-2\deg \widetilde h_d +m_d+4\deg S_d=-6d(d-2)+m_d+12(d-1)^2.
\end{equation}

On the other hand, the curve $I_L\subset L\times \mathbb P^2$ is a complete intersection of $H_{C_L}$ and the smooth surface $C_L$.
The Picard group $Pic(L\times \mathbb P^2)$ is a free abelian group generated by divisirs $A=pr_1^{-1}(t)$ and $B=pr_2^{-1}(\mathbb P^1)$, where $t$ is a point in $L$ and $\mathbb P^1$ is a line in $\mathbb P^2$. The canonical class $K_{L\times \mathbb P^2}=-2A-3B$, $C_L\in|A+dB|$, and $H_{C_L}\in |3A+3(d-2)B|$. Let us restrict the divisors $A$ and $B$ to the surface $C_L$. It follows from the adjunction formula that
$K_{C_L}=-A+(d-3)B$. Besides, we have $I_L\in |3A+3(d-2)B|$,
$(A,A)_{C_L}=(A,A,A+dB)_{L\times\mathbb P^2}=0, \quad (A,B)_{C_L}=(A,B,A+dB)_{L\times\mathbb P^2}=d$,
and $(B,B)_{C_L}=(B,B,A+dB)_{L\times\mathbb P^2}=1$.

Therefore we have
$ 2(g_d-1)=  (I_L,I_L+K_{C_L})_{C_L}-2\deg \mathcal S_d=  6(4d^2-13d+8)$.
Сomparing this equality with (\ref{g1}), we obtain that
$m_d=6(d-3)(3d-2)$. \qed \\

Let $L$ be a line in $\mathbb P^{K_d}$, $L\not\subset \mathcal S_d$, and $I_L=h_d^{-1}(L)$. We call $H_L=pr_2(I_L)\subset \mathbb P^2$ the {\it Hesse curve of the pencil} $\{ C_{\overline a}\}_{\overline a\in L}$.

Since for a generic Lefschetz pencil $\{ C_{\overline a}\}_{\overline a\in L}$ the Hesse curve $H_L$ is irreducible, simple calculations show that the number of elements of a generic Lefschetz pencil having an inflection point at a base point of the pencil is less than or equal to three.  Besides, we have  $\deg H_L=(H_L,\mathbb P^1)_{\mathbb P^2}=(I_L,B)_{C_L}=6(d-1)$. Therefore we have

\begin{thm} The Hesse curve $H_L$ of a generic Lefschetz pencil $\{ C_{\overline a}\}_{\overline a\in L}$ of plane curves of degree
$d\geq 3$ has the following properties:
\begin{itemize} \item[$(i)$]
$\deg H_L=6(d-1)$ and its genus is equal to $3(4d^2-13d+8)+1$;
\item[$(ii)$]
 $H_L$  has $d^2$ singular points of multiplicity three at the base points of the pencil $\{ C_{\overline a}\}_{\overline a\in L}$
and $3(d-1)^2$ ordinary nodes at the singular points of the degenerate fibres of the pencil;
\item[$(iii)$]
 for each $\overline a\in L$ the intersection $H_L\cap C_{\overline a}$ consists of the inflection points of $C_{\overline a}$ and the base points of the pencil, and if $p\in H_L\cap C_{\overline a}$ is the $2$-tuple inflection point of $C_{\overline a}$, then $H_L$ and $C_{\overline a}$ touch each other at $p$.
\end{itemize}
 \end{thm}

\ifx\undefined\bysame
\newcommand{\bysame}{\leavevmode\hbox to3em{\hrulefill}\,}
\fi

\ifx\undefined\bysame
\newcommand{\bysame}{\leavevmode\hbox to3em{\hrulefill}\,}
\fi

\end{document}